\documentclass{article}
\usepackage{amssymb,amsmath}

\newtheorem{theorem}{Theorem}[section]

\newtheorem{hypothesis}{Conjecture}[section]

\begin{document}

\title{Progress in number theory in the years 1998-2009}

\author{Adam Grygiel\footnote{The author was partially supported by Dean of the Faculty of Mathematics and Computer Science of the University of \L \'od\'z, and the European Social Fund and Budget of State implemented under the Integrated Regional Operational Programme (Project: GRRI-D).
}}

\date{}

\maketitle

\begin{abstract}
We summarize the major results in number theory of the last decade.
\end{abstract}

\section{Introduction}
The purpose of the present paper, originally published in Polish (see \cite{wm1}), is to review briefly nine spectacular achievements belonging 
to the theory of numbers from the years 1998-2009. We classify 
these results in the following subjects according to \textit{Mathematical 
Reviews}: 
\begin{itemize}
\item Elementary number theory; 
\item Sequences and sets of integers; 
\item Diophantine equations; 
\item Analytic number theory;  
\item Computational number theory.
\end{itemize}
I would like to thank the anonymous referee for his remarks improving the paper, and Professors Jerzy Browkin and Andrzej Schinzel for their valuable comments and advice. I am also grateful to Professor Kevin Ford for kindly correcting the statement of his result in the second section, about the Carmichael Conjecture. 

\section{Elementary number theory}

For \textit{Euler's totient function} $\phi$ the following formula holds
\[
\phi (n)=
n\prod_{i=1}^{k}\left( 1-\frac{1}{p_i}\right), 
\]
in which $p_1,\ldots , p_k$ denote distinct prime factors of $n$. (Let us remember that by $\phi (a)$ we denote the number of positive integers less than or equal to $a$ that are coprime to $a$.) 
In 1922, R.D. Carmichael formulated in \cite{c} an hypothesis, which asserts that there is no integer $m$ such that the equation $\phi (x)=m$ has exactly one solution; this conjecture is still unsubstantiated. Denote by $\mathcal{V}$ the set of values taken by Euler's totient function.
In 1998, K. Ford proved in \cite{f2} that a prospective counter-example abolishing Carmichael's conjecture satisfies $m>10^{10^{10}}$, and that if there exists such counter-example, then the set $A$ of these counter-examples has positive relative lower density, that is,
\[
\liminf_{N \to \infty}\frac{|A\cap [1,N]|}{|\mathcal{V}\cap [1,N]|}>0.
\]

W. Sierpi\'nski has formulated the following hypothesis (see \cite[Ch. VI, 1, p. 252]{6}).
\begin{hypothesis}[Sierpi\'nski 1950]\label{hs1}
For any integer $s\geq 2$, there exists such an integer 
$m$ that the equation $\phi (x)=m$ has exactly $s$ 
solutions.
\end{hypothesis}

K. Ford and S. Konyagin have proved in  \cite{fk1} the following 
theorem.

\begin{theorem}[Ford, Konyagin 1999]\label{tfk}
Sierpi\'nski's Conjecture on Euler's totient function holds for 
even integers.
\end{theorem}

Soon after this breakage, K. Ford  has proved in \cite{24} Conjecture \ref{hs1}. 
He has used many deep results from sieve 
theory, e.g. Chen's theorem from \cite{25}, which asserts that 
there exist infinitely many prime numbers $p$ such that $p+2$ 
is a product of at least two prime numbers. 

\section{Sequences and sets of integers}

P. Erd\H os (see \cite[p. 11]{erd}) has offered \$3000 for a solution of the following 
still unsolved problem.
\begin{hypothesis}[Erd\H os 1980]\label{et}
Any subset $A$ of the positive integers such that 
$\sum_{n\in A}\frac{1}{n}=\infty$ contains an arithmetic 
progression of length $k$, for all $k$.
\end{hypothesis}

It is known due to Euler that the series of the inverses 
of all prime numbers divergents. In 1939, J.G. van der 
Corput  proved in \cite{12} that there exist infinitely 
many arithmetic progressions of prime numbers of length 
$3$. In 1975, E. Szemer\'edi proved in  \cite{10}, in a 
combinatorial way, that any subset of positive integers, with
positive upper density, contains 
arithmetic progression of length $k$, for all $k$. 
Unfortunately, the set $\mathcal{P}$ of all prime numbers 
has upper density zero.  W.T. Gowers  has extended Szemer\'edi's theorem, 
by using Fourier analysis, to the following result from \cite{go}.
\begin{theorem}[Gowers 2001, Fields Medal 1998]
\label{tgo}
The maximal length $r_l(n)$ of a progression of integers 
not exceeding $n$, containing no arithmetic progression 
of length $l$, satisfies
\[
r_l(n)=O\left( \frac{n}{(\ln \ln n)^{c_l}}\right) , 
\quad c_l=2^{-2^{l+9}}.
\]
\end{theorem}

In 2005, B. Green  proved in \cite{77}, by the same method, 
that any set $A\subset\mathcal{P}$ with positive 
relative upper density, that is, satisfying the condition
\[
\limsup_{N \to \infty}\frac{|A\cap [1,N]|}{|\mathcal{P}
\cap [1,N]|}>0,
\]
 contains infinitely many arithmetic 
progressions of length $3$. By using ergodic methods, jointly with T. Tao he has generalized this fact in \cite{5} to arbitrarily long arithmetic progressions, solving 
Conjecture \ref{et} in the crucial case, when $A=\mathcal{P}$.

\begin{theorem}[Green, Tao 2008, Fields Medal 2006]
\label{gt}
Any set $A\subset\mathcal{P}$ with positive relative upper 
density contains infinitely many arithmetic progressions of 
length $k$, for all $k$.
\end{theorem}

Let $N\in \mathcal{P}$ be a sufficiently large. Let us denote
\[
W=\prod_{\substack{p \in \mathcal{P}\\ p\leq \ln \ln N}}p.
\]
Let us remember that \textit{von Mangoldt's function} $\Lambda$ is given by
\begin{displaymath}
\Lambda (n)=
\begin{cases}
\ln p & \textup{if} \ n=p^l  \ \textup{for some} \ p 
\in \mathcal{P} \ \textup{and positive integer} \  l,\\
0 & \textup{in the opposite case},
\end{cases}
\end{displaymath}
and consider the following modification of this function:
\begin{displaymath}
\tilde \Lambda (n)=
\begin{cases}
\frac{\phi (W)}{W}\ln (Wn+1) & \textup{if} \ Wn+1 \in 
\mathcal{P},\\
0 & \textup{in the opposite case}.
\end{cases}
\end{displaymath}
A key point in the proof of Theorem \ref{gt} is a lower 
evaluation of the expression
\[
\frac{1}{N^2}\sum_{n=1}^{N}\sum_{r=1}^{N}\tilde \Lambda 
(n) \tilde \Lambda (n+r)\cdots \tilde \Lambda (n+(k-1)r).
\]
This evaluation implies that there 
exists in $\mathcal{P}$ an arithmetic progression of the form
\[
Wn+1,W(n+r)+1,\ldots , W(n+(k-1)r)+1.
\] 

The proof of Theorem \ref{gt} does not follow as to design 
arithmetic progressions in $\mathcal{P}$ of a given length. 
In 2008, J. Wr\'oblewski and R. Chermoni found the longest 
currently known such progression:
\[
6171054912832631 + 366384 \times  223092870 \times  n, 
\quad n = 0,\ldots , 24.
\]

T. Tao and T. Ziegler have proved in \cite{11}, by using ergodic 
theory, the following generalization of Theorem \ref{gt}.
\begin{theorem}[Tao, Ziegler 2008]\label{tz}
If $P_1,\ldots,P_k \in \mathbb{Z}[m]$ are integer-valued 
polynomials such that $P_1(0)=\cdots=P_k(0)=0$, then 
any subset of $\mathcal{P}$ with positive relative
upper density contains infinitely many sequences of 
the form  $n+P_1(m),\ldots ,n+P_k(m)$, with $m>0$.
\end{theorem}

\section{Diophantine equations}

E. Catalan has formulated in \cite{0}  the following hypothesis.
\begin{hypothesis}[Catalan 1844]\label{hc}
Catalan's equation 
\[
x^p-y^q=1,
\]
has no  solutions in integers $x,y,p,q>1$ other than $3^2-2^3=1$.
\end{hypothesis}

The case of $q=2$ of Conjecture \ref{hc} was solved  in \cite{l} 
by V.A. Lebesgue in 1850. In 1964, Chao Ko proved in  \cite{k} 
Conjecture \ref{hc} for $p=2$. In 1976, R. Tijdeman 
proved in \cite{ti}, by using Baker's method of estimates for linear forms 
of logarithms, that Catalan's equation has only finitely many 
solutions. These results were clearly presented in
\cite{pr} by P. Ribenboim.

In 1990, K. Inkeri (see \cite{i1,i2}) proved the following result 
(the so called \textit{Inkeri's criterion}). Let $p,q\in \mathcal{P}$ 
be odd integers. If Catalan's equation has a solution in 
integers $x,y>1$, then the following alternative holds: 
$p^{q-1}\equiv 1$ (mod $q^2$) or $q$ divides the class number 
of a number field $L$ defined as follows:
\[
L=
\begin{cases}
\mathbb{Q}\left( \sqrt{-p}\right) & \textup{if} \ p\equiv 
3 \ (\textup{mod} \ 4),\\
\mathbb{Q}\left( e^{2i\pi /p}\right) & \textup{in the opposite 
case}.
\end{cases}
\]

In 2003, P. Mih$\check{\textup{a}}$ilescu  proved in \cite{mi}
that the second term of the alternative in Inkeri's 
criterion one can drop. More precisely, he has proved the following 
theorem.
\begin{theorem}[Mih$\check{\textup{a}}$ilescu 2003]\label{tmi}
Let $p,q\in \mathcal{P}$ be odd integers. If Catalan's 
equation has a solution in integers $x,y>1$, then 
$p^{q-1}\equiv 1$ (mod $q^2$) and $q^{p-1}\equiv 1$ (mod $p^2$).
\end{theorem}
 
A pair of odd integers $p,q\in \mathcal{P}$, satisfying 
both congruences in Theorem \ref{tmi}, is called a
\textit{double Wieferich pair}. There are currently only six 
such pairs known. In 2004, P. Mih$\check{\textup{a}}$ilescu 
 proved in \cite{99} Conjecture \ref{hc}. A crucial role in his 
 proof is played by the condition 
 $p\not \equiv 1$ (mod $q$), 
which he has proved, by using the 
double Wieferich pair condition. The original proof was much 
improved by Y. Bilu (see \cite{yb1, yb2}).

\section{Analytic number theory}

A. Schinzel has formulated in \cite{9} the following general 
hypothesis, known as \textit{Schinzel's Hypothesis H}.
\begin{hypothesis}[Schinzel 1958]\label{sch}
Let $f_1,\ldots,f_k\in \mathbb{Z}[m]$ be irreducible, 
integer-valued polynomials, with positive leading 
coefficients. If for every $q \in \mathcal{P}$ we have 
$q\nmid f_1(m)\cdots f_k(m)$ for some $m \in \mathbb{Z}$, 
then $f_1(n),\ldots,f_k(n)\in \mathcal{P}$ for infinitely 
many positive integers $n$. 
\end{hypothesis}
In 1961, A. Schinzel proved in \cite{21} that Conjecture \ref{sch} 
implies Conjecture \ref{hs1}. Jointly with W. Sierpi\'nski he has 
also deduced in  \cite{9} many other interesting 
corollaries from Conjecture \ref{sch}, e.g. that there exist arbitrarily long arithmetic 
progressions of consecutive prime numbers (which is still an open 
problem). The longest currently known such progression, of length 10, 
was found by a group associated with M. Toplic in 1998.

Classical Dirichlet's theorem on prime numbers in 
arithmetic progressions says that, if $f(m)=bm+a$, where 
$a,b\in \mathbb{Z}$, $a\not =0$, $b\geq 1$, and $(a,b)=1$, 
then $f(n)\in \mathcal{P}$ for infinitely many integers $n$. In 1978, H. Iwaniec 
proved in \cite{7a} that $n^2+1$ is a product of at least two prime numbers for 
infinitely many integers $n$. 
We know currently no polynomial of degree greater than $1$, 
in one variable, which would represent infinitely many prime 
numbers. Also for $k>1$  Conjecture \ref{sch} is 
completely open problem, even for linear polynomials. 

It is known due to Euler (which is the statement of Fermat's theorem on sums of two 
squares) that a prime number $p>2$ is a sum of two squares 
of integers, if and only if $p\equiv 1 \ (\textup{mod} \ 4)$. 
In particular, $m^2+n^2\in \mathcal{P}$ for infinitely many 
integers $m,n$. In 1974, H. Iwaniec  generalized in \cite{7} the 
last fact to polynomials of degree $2$, in two variables, 
satisfying some natural assumptions. In 1997, E. Fouvry 
and H. Iwaniec proved in \cite{fi}, by using sieve methods, 
that $m, m^2+n^2\in \mathcal{P}$ for infinitely many integers 
$m,n$. J. Friedlander and H. Iwaniec have proved in  \cite{88}, by 
using Bombieri's asymptotic sieve, that
\[
\underset{m^2+n^4\leq x}{\sum\sum}\Lambda (m^2 + n^4)
\sim c x^{3/4}, \quad 
c=\frac{4}{\pi}\int_0^1(1-t^4)^{1/2}\textup{dt}.
\]
This implies the following theorem.
\begin{theorem}[Friedlander, Iwaniec 1998] 
$m^2 + n^4\in \mathcal{P}$ for infinitely many integers $m,n$.
\end{theorem}

D.R. Heath-Brown has proved in \cite{66}, by the same method, the 
following theorem.
\begin{theorem}[Heath-Brown 2001]
$m^3 + 2n^3\in \mathcal{P}$ for infinitely many integers 
$m,n$.
\end{theorem}

Let $p_n$ be the $n$-th prime number. Let us denote
\[
\triangle _1=
\liminf _{n\to \infty}\frac{p_{n+1}-p_n}{\ln p_n}.
\]
Until quite lately, the best known evaluation was 
$\triangle _1\leq 0.2484$, proven in \cite{ma} by H. Maier 
 in 1988. D.A. Goldston, J. Pintz and 
C.Y. Y\i ld\i r\i m \cite{gpy} have reached a great breakage in 
this area. They have proved in \cite{gpy}, by using Selberg's sieve methods, 
the following theorem.
\begin{theorem}[Goldston, Pintz, Y\i ld\i r\i m 2009]
\label{tgpy} $\triangle _1=0$.
\end{theorem}
 
D.A. Goldston, J. Pintz and C.Y. Y\i ld\i r\i m have also proved in 
\cite{gpy1}  the following reinforcement of Theorem \ref{tgpy},
\[
\liminf _{n\to \infty}\frac{p_{n+1}-p_n}{(\ln p_n)^{1/2
}
(\ln \ln p_n)^2}<\infty .
\]
 
For fixed $\varepsilon >0$ let us take sufficiently large 
integers $N,k$, and define $h=\varepsilon \ln N$. Let us denote 
$\mathcal{H}_k=\{h_1,\ldots ,h_k\}$, where $h_1,\ldots ,h_k$ 
are integers and $1\leq h_1<\ldots <h_k\leq h$, and consider, 
when a polynomial $P_{\mathcal{H}_k}$ given by 
$P_{\mathcal{H}_k}(n)=(n+h_1)\cdots (n+h_k)$ has $k+l$ or less distinct prime factors, 
where $0\leq l\leq k$. Set
\[
\Lambda _R(n;\mathcal{H}_k,l)=\frac{1}{(k+l)!}\sum_{
\substack{d|P_{\mathcal{H}_k}(n) 
\\ d\leq R}}\mu (d) \left(\ln \frac{R}{d}\right)^{k+l},
\]
where $R$ is a suitable real parameter; let us remember that 
\begin{displaymath}
\mu (n)=
\begin{cases}
(-1)^k & \textup{if} \ n \ \textup{is a product of $k$ 
distinct prime numbers},\\
0 & \textup{for other} \ n>1.
\end{cases}
\end{displaymath}
Let us consider the following modification of von 
Mangoldt's function $\Lambda$:
\begin{displaymath}
\hat \Lambda (n)=
\begin{cases}
\ln n & \textup{if} \ n \in \mathcal{P},\\
0 & \textup{in the opposite case}.
\end{cases}
\end{displaymath}
A key point in the proof of Theorem \ref{tgpy} 
is a lower evaluation of the expression
\[
\sum _{n=N+1}^{2N}\left( \sum _{1\leq h_0\leq h}\hat 
\Lambda (n+h_0)-\ln (3N)\right) \sum_{1\leq h_1< \ldots 
< h_k\leq h}\Lambda _R(n;\mathcal{H}_k,l)^2.
\]
This evaluation implies that the 
interval $(n,n+h)$ contains at least two prime numbers 
for infinitely many integers $n$.

\section{Computational number theory}

In 1983, L.M. Adleman, C. Pomerance and R.S. Rumely 
 published in \cite{a1} an algorithm for determining whether a 
given integer $n>1$ is a prime. The time complexity of 
their algorithm equals
\[
(\ln n)^{O\left( \ln \ln \ln n\right)}.
\]

In 2004, M. Agrawal, N. Kayal and N. Saxena
published in  \cite{a}  the first polynomial-time primality test. Their  
algorithm (the so called \textit{AKS primality test}) executes the 
following steps.

\begin{enumerate}
\item If $n=a^b$ for some integers $a,b>1$, output 
,,composite''.
\item Find the smallest $r$ such that $o_r(n)>
(\log _2 n)^2$, where $o_r(n)$ denotes the smallest positive integer $k$ such that  $n^k\equiv 1$ 
(mod $r$), and $\log _2$ means a logarithm to the base $2$.
\item If $1<(a,n)<n$ for some $a\leq r$, output 
,,composite''.
\item If $n\leq r$, output ,,primes''.
\item For $1\leq a\leq \lfloor \sqrt{\phi (r)}\log _2 n 
\rfloor$ do: $\textup{if} \ (X+a)^n\not =X^n+a \
\textup{in} \ \left(\mathbb{Z}/n\mathbb{Z}\right)[X]/
(X^r-1), \ \textup{output ,,composite''}$.
\item Output ,,primes''.
\end{enumerate}

\begin{theorem}[Agrawal, Kayal, Saxena 2004]\label{aks}
The AKS primality test returns ,,primes'' if and only if 
$n\in \mathcal{P}$. The time complexity of this algorithm 
equals $O\left((\ln n)^{\frac{21}{2}+\varepsilon }\right)$, for any 
$\varepsilon >0$.
\end{theorem}

The main difficulty in the proof of correctness of the AKS primality test lies in the implication that, if the above algorithm returns ,,primes'', then $n\in \mathcal{P}$. It was shown elementarily (see \cite[Lemma 4.3]{a}) that there exists $r\leq \max \{3,\lceil (\log_2 n)^{5} \rceil \}$ such that $o_r(n)>(\log_2 n)^2$. Since $o_r(n)>1$, there exists such prime factor $p$ of $n$ that $o_r(p)>1$. A further part of the proof is based on the equality 
\[
x^r-1=\prod _{d|r}\Phi _d(x),
\]
in which $\Phi _d\in \left(\mathbb{Z}/p\mathbb{Z}\right)[x]$ denotes the $d$-th cyclotomic polynomial. Let  $h\in \left(\mathbb{Z}/p\mathbb{Z}\right)[x]$ be an irreducible factor of $\Phi _r$. Then, the ring $\mathbb{F}:=\left(\mathbb{Z}/p\mathbb{Z}\right)[x]/(h(x))$ is a finite field of order $p^d$, where $d$ is the degree of $h$. A key point in the proof of Theorem \ref{aks} is a lower and upper evaluation of the order of the cyclic subgroup $\mathbb{F}$ generated multiplicatively by the elements
\[
x,x+1,x+2,\ldots, 
x+ \lfloor \sqrt{\phi (r)}\log _2 n \rfloor , 
\]
under the assumption that $n$ is not a power of $p$. It follows from these evaluations  that $n=p$. 

H.W. Lenstra, Jr. and C. Pomerance have
modificated in \cite{lp} the AKS primality test  for obtaining a deterministic primality test with the time 
complexity $O\left((\ln n)^{6+\varepsilon }\right)$, for 
any $\varepsilon >0$. 

\bibliographystyle{amsplain}

\vskip10pt

\textit{Faculty of Mathematics and Computer Science, University of 
\L\'od\'z, Banacha 22, 90-238 \L\'od\'z, Poland}\\
A.Grygiel@math.uni.lodz.pl

\end{document}